\documentclass[review]{article}

\usepackage{lineno,hyperref}
\modulolinenumbers[5]

\usepackage{amssymb}
\usepackage{amsthm}
\usepackage{amsmath}

\newcommand{\SA}{\mathcal{A}}
\newcommand{\SD}{\mathcal{D}}
\newcommand{\SL}{\mathcal{L}}

\theoremstyle{plain}
\newtheorem{theorem}{Theorem}[section]

\theoremstyle{definition}

\theoremstyle{remark}
\newtheorem{remark}{Remark}
\newtheorem{notation}{Notation}

\begin{document}


\title{Generalized Pseudospectral Method and Zeros of Orthogonal Polynomials}



\date{}
\author{Oksana Bihun\footnote{Corresponding author, obihun@uccs.edu}, Clark Mourning\\
\small{Department of Mathematics,
University of Colorado, Colorado Springs,}\\
\small{1420 Austin Bluffs Pkwy,
Colorado Springs, CO 80918, USA}}

\maketitle

\begin{abstract}
Via a generalization of the pseudospectral method for numerical solution of differential equations,  a family of nonlinear algebraic identities satisfied by the zeros of  a wide class of orthogonal polynomials is derived. The generalization is based on a modification of pseudospectral matrix representations of linear differential operators proposed in the paper, which allows these representations to depend on two, rather than one, sets of interpolation nodes.
The identities hold for every polynomial family $\{p_\nu(x)\}_{\nu=0}^\infty$   orthogonal with respect to a measure supported on the real line that satisfies some standard assumptions, as long as the polynomials in the family satisfy differential equations $\SA p_\nu(x) =q_\nu(x) p_\nu(x)$, where $\SA$ is a linear differential operator and each $q_\nu(x)$ is a polynomial of degree at most $n_0 \in \mathbb{N}$; $n_0$ does not depend on $\nu$.  The proposed identities generalize known identities for classical and Krall orthogonal polynomials, to the case of the nonclassical orthogonal polynomials that belong to the class described above. The generalized pseudospectral representations of the differential operator $\SA$ for the case of the Sonin-Markov orthogonal polynomials, also known as generalized Hermite polynomials, are presented. The general result is illustrated by new algebraic relations satisfied by the zeros  of the Sonin-Markov polynomials.
\end{abstract}

\textbf{Keywords:}
 Pseudospectral methods; Spectral methods; Zeros of polynomials; Nonclassical orthogonal polynomials; Sonin-Markov polynomials; Generalized Hermite polynomials.
 
\textbf{MSC2010} 33C45; 65L60; 33C47; 26C10.



\section{Introduction and Main Results}
\label{Sec:Intro}

\subsection{Summary of results}
In this paper we identify a class of  algebraic relations satisfied by the zeros of  a wide class of orthogonal polynomials.  To prove the identities, we generalize the notion of pseudospectral matrix representations of linear differential operators, by allowing these representations to depend on two, rather than one, sets of interpolation nodes.
 The identities hold for all polynomials $\{p_\nu(x)\}_{\nu=0}^\infty$ orthogonal with respect to a measure satisfying some standard assumptions, as long as they satisfy the differential equations 
\begin{equation}
\SA p_\nu(x) =q_\nu(x) p_\nu(x),
\label{EigenFunctionEquation}
\end{equation}
 where $\SA$ is a linear differential operator and each $q_\nu(x)$ is a polynomial of degree at most $n_0 \in \mathbb{N}$; $n_0$ does not depend on $\nu$. This includes classical orthogonal polynomials~\cite{Szego, NikiforovUvarov}, polynomials in the Askey scheme~\cite{KoekSwart} and Krall polynomials~\cite{AMKrallBook02}, as well as additional classes of nonclassical orthogonal polynomials. If applied to classical or Krall orthogonal polynomials, the proposed result reduces to known identities~\cite{Ahmed,  Sasaki15, Bihun2016}, see Subsection~\ref{MainRes}. 


The motivation of this study stems from understanding that zeros of
orthogonal polynomials play an important role in mathematical physics, numerical analysis and related areas. For example, zeros of some orthogonal polynomials are equilibria of important N-body
problems \cite{CalogeroBook, C2012, CBookIsoch, BPSbook2011}. They transpire as building blocks of remarkable isospectral matrices \cite{BCGenHyperg2014, BCAskey2014, BCGenBHyperg2015, BCqAskey2016} and play an important role in construction of highly accurate approximation schemes for numerical integration \cite{ MastMilo, MastronardiOcco00, MastOcco04, OccoRusso2014}.

To prove algebraic identities satisfied by the zeros of the polynomials $\{p_\nu(x)\}_{\nu=0}^\infty$, we generalize and relate the notions of spectral and pseudospectral matrix representations of linear differential operators used in the corresponding numerical methods for solving differential equations.
The standard pseudospectral matrix representations of linear differential operators are based on Lagrange collocation on the real line. These representations were proposed
by Calogero  in the context of numerical solving of eigenvalue and boundary value problems for linear ODEs~\cite{Calogero82, Calogero83}, see also~\cite{CF84, BCP90, Durand83, CalogeroBook, SpectralMethodsBook}. Mitropolsky, Prykarpatsky and Samoylenko set up a general algebraic-projection framework based on Calogero's method
and furthered its applications to solution of evolution equations in Mathematical Physics~\cite{MPS88}. The convergence analysis of  Calogero's method was studied in~\cite{CEJMMirekCo}.

The standard pseudospectral method was utilized in~\cite{Bihun2016} to prove new properties of the zeros Krall polynomials. While Krall polynomials are eigenfunctions of linear differential operators, the polynomial families considered in this paper  satisfy differential equations~\eqref{EigenFunctionEquation}  with $q_n(x)$ being \textit{polynomials (as opposed to eigenvalues)} of degree $n_0$ that does not depend on $n$. To prove new properties of the zeros of the latter polynomial families, we propose a \textit{generalization} of the standard pseudospectral method.

Because, in general,  the differential operator $\SA$ in~\eqref{EigenFunctionEquation} raises the degree of polynomials by a summand of $n_0$, the standard pseudospectral method does not allow \textit{exact} discretization of differential equations~\eqref{EigenFunctionEquation}.  
The main idea of the proposed generalization is to construct Lagrange collocation type matrices that  \textit{exactly}
represent linear differential operators acting between spaces of polynomials of different degrees. 
In this paper, the last goal is achieved by allowing the matrix representations to depend on two rather than one set of interpolation nodes. We thus find  \textit{exact} discretizations of the differential equations~\eqref{EigenFunctionEquation} satisfied by the polynomials $\{p_\nu(x)\}_{\nu=0}^\infty$; the discretizations are constructed using the zeros of these polynomials as the nodes. By comparing the generalized pseudospectral and the generalized spectral matrix representations of the differential operator $\SA$, we derive a family of algebraic identities satisfied by the zeros of polynomials $\{p_\nu(x)\}_{\nu=0}^\infty$.

The proposed generalization of the pseudospectral matrix representations of linear differential operators has applications beyond those outlined in this paper; one such application allows to simplify the process of incorporation of initial or boundary conditions into linear systems that discretize certain ODEs, see the discussion in Section~\ref{Sec:Outlook}.

We illustrate the general result of the main Theorem~\ref{PropZerosThmG} by applying it to the case of the Sonin-Markov polynomials, also known as generalized Hermite polynomials, see~~\cite{Kis, Chihara, MastMilo} and references therein. These zeros play an important role in 
the computation of integrals of singular or oscillatory functions~\cite{MastronardiOcco00, MastOcco04} as well as extended Lagrange interpolation on the real line~\cite{OccoRusso2014}.

To prove the identities satisfied by the zeros of the Sonin-Markov polynomials stated in Theorem~\ref{PropZerosThmSM}, we compute the generalized  pseudospectral matrix representations of the differential operators $d/dx$ and $d^2/dx^2$ as well as the differential operator $\SD$ associated with the Sonin-Markov family, see~\eqref{DSoninMarkov} of~\eqref{diffEqSoninMarkov} and Subsection~\ref{sec:DcSM}, assuming that the interpolation nodes are   zeros of  the Sonin-Markov polynomials. The formulas for these matrix representations as well as the identities of Theorem~\ref{PropZerosThmSM} have been verified using programming environment MATLAB, for several particular values of the relevant paramenets (the degree $N$ and the parameter $\beta$, see~Theorem~\ref{PropZerosThmSM} and definitions~\eqref{SMdf},~\eqref{alpha}).

\subsection{The orthogonal polynomial family $\{p_{\nu}(x)\}_{\nu=0}^\infty$}

Let $\{p_\nu(x)\}_{\nu=0}^\infty$ be a sequence of polynomials orthogonal  with respect to a measure $\omega$ and the corresponding inner product $\langle f,g\rangle=\int fg \;d \omega$. We denote the norm associated with this inner product by $\|\cdot \|$, that is, $\|f\|^2=\int f^2 \; d\omega$.  Assume that $\omega$ is a Borel measure with support on the real line satisfying the following three conditions:
\begin{itemize}
\item[(a)] $\omega$ is positive;
\item[(b)] all its moments $\int x^\nu  \,d\omega$ exist and are finite;
\item[(c)] $\omega$ has infinitely many points in its support $I=\operatorname{supp} \omega$.
\end{itemize}
Under the above assumptions on the measure $\omega$, the zeros of each polynomial $p_\nu,\; \nu\geq 1,$ are real, simple and belong to the convex hull of the support of $\omega$, see for example~\cite{MastMilo}.
 
 \begin{notation}
 \label{Notation1}
Here and throughout the rest of the paper $N$  denotes a fixed integer strictly larger than 1, while $n_0$ is a fixed nonnegative integer. The small Greek letter $\rho$ denotes an index that may take values $N$ or $N+n_0$. The small Greek letter $\nu$ denotes an integer index that usually takes values $0,1, 2 \ldots$,  unless otherwise indicated. The small Latin letters $n,m,j,k$ etc. denote integer indices that usually run from $1$ to $N$ or from  $1$ to $N+n_0$, see~\eqref{DiffOperatorA} and~\eqref{EigenFunctionEquation}, thus we indicate the range of the indices each time they are used. We reserve  the letter  $\ell$  to denote polynomials in Lagrange interpolation bases.
\end{notation}
 
 Let $\mathbb{P}^\nu$ denote the space of all algebraic polynomials with real coefficients of degree at most $\nu$.  Assume that for each $\nu$, the polynomials $\{p_j(x)\}_{j=0}^\nu$ form a basis of $\mathbb{P}^\nu$.
 Let $\SA$ be a linear differential operator acting on functions of one variable. Assume that $\SA$ has the property
 \begin{equation}
 \SA \mathbb{P}^\nu \subseteq \mathbb{P}^{\nu+n_0}
 \label{DiffOperatorA}
 \end{equation} for all $\nu$. Recall that $n_0$ is a fixed nonnegative integer; it does not depend on $\nu$. For example, the differential operator $\mathcal{D}=a_0+a_1(x) \frac{d}{dx}+\ldots+a_q(x)\frac{d^q}{dx^q}$ with $q \in \mathbb{N}$ and $a_j(x) \in \mathbb{P}^j$ for all $j=1,2, \ldots, q$ has property~\eqref{DiffOperatorA} with $n_0=0$, while the operator $x^2 \mathcal{D}$ has property~\eqref{DiffOperatorA} with $n_0=2$.
 
Suppose that the orthogonal polynomials $\{p_\nu(x)\}_{\nu=0}^\infty$ satisfy differential equations~\eqref{EigenFunctionEquation}.

 \subsection{Generalized pseudospectral and spectral matrix representations of linear differential operators}
 
In this Subsection we  introduce the notions of a generalized pseudospectral and a generalized spectral matrix representations of the linear differential operator  operator $\SA$ introduced in the previous Subsection.  Note that, in general,  the definitions of these generalized matrix representations hold for any linear differential operator $\SD$, which may or may not satisfy property~\eqref{DiffOperatorA}, and for $n_0$ being a (positive or negative) integer such that $N+n_0\geq 1$.
 
The definition of the standard pseudospectral $N \times N$ matrix representation $\tilde{A}^c$ of $\SA$, see~\cite{Bihun2011,BihunPrytula2010,Bihun2016,SpectralMethodsBook}, is motivated by a search for an exact discretization of a differential equation 
\begin{equation}
\SA u=f,
\label{Au=f}
\end{equation}
under the assumption that it possesses a polynomial solution $u \in \mathbb{P}^{N-1}$. More precisely, choose a vector $\vec{x}^{(N)}=\left( x^{(N)}_1, \ldots, x^{(N)}_N\right)$ of distinct real nodes and define the isomorphism $\pi_N: \mathbb{P}^{N-1} \to \mathbb{R}^N$ by 
\begin{equation}
\pi_N g=\left( g\left(x^{(N)}_1\right), \ldots, g\left(x^{(N)}_N\right)\right).
\label{piN}
\end{equation}
The inverse of $\pi_N$ is of course given in terms of the standard Lagrange interpolation basis $\{\ell_{N-1, j}(x)\}_{j=1}^N$ of $\mathbb{P}^{N-1}$ constructed using the nodes $\vec{x}^{(N)}$: for every vector $\vec{g}^{(N)}=\left( g^{(N)}_1, \ldots, g^{(N)}_N\right) \in \mathbb{R}^N$,
$$
\pi_{N}^{-1} \vec{g}^{(N)}=\sum_{j=1}^N g^{(N)}_j \ell_{N-1,j}(x).
$$
 The standard pseudospectral $N \times N$ matrix representation $\tilde{A}^c$ of the differential operator $\SA$ is defined as the unique $N \times N$ matrix that satisfies the condition $$\pi_N \SA g= \tilde{A}^c \pi_N g$$ for all $g \in \mathbb{P}^{N-1}$. Note that the superscript `$c$' in the notation of the matrix $\tilde{A}^c$ stands for ``collocation''  in the spectral collocation method for numerical solving of differential equations, also known as the pseudospectral method~\cite{SpectralMethodsBook}.  It is not difficult to conclude that the matrix $\tilde{A}^c$ is given componentwise by 
 $$\tilde{A}^c_{kj}=\SA \ell_{N-1,j}(x)\big|_{x=x_k^{(N)}},$$ 
 where $1\leq k,j \leq N$, see~\cite{Calogero82, Calogero83, MPS88, CEJMMirekCo, BihunPrytula2010, Bihun2011, Bihun2016}.
The last definition implies that a vector $\vec{u}^{(N)} \in \mathbb{R}^N$ solves the linear system 
\begin{equation}
\tilde{A}^c \vec{u}^{(N)}=\pi_N f
\label{Au=fApprox}
\end{equation}
if and only if the polynomial $u(x)=(\pi_N)^{-1}\vec{u}^{(N)}=\sum_{j=1}^N {u}^{(N)}_j \ell_{N-1, j}(x)$ solves  ODE~\eqref{Au=f}. Of course, if ODE~\eqref{Au=f} does not possess a polynomial solution $u \in \mathbb{P}^{N-1}$,  linear system~\eqref{Au=fApprox} \textit{approximates} ODE~\eqref{Au=f} and its solution $\vec{u}^{(N)}$, if it exists, allows to construct an \textit{approximate} solution of ODE~\eqref{Au=f} given by $u(x)=(\pi_N)^{-1}\vec{u}^{(N)}$.

Let us now generalize this notion of pseudospectral matrix representation of the differential operator $\SA$ to take advantage of its property~\eqref{DiffOperatorA}.

In addition to the interpolation nodes $\vec{x}^{(N)}$, consider another vector of distinct real nodes  $\vec{x}^{(N+n_0)}=\left( x^{(N+n_0)}_1, \ldots, x^{(N+n_0)}_{N+n_0}\right)$. In short, we will work with two vectors of nodes $\vec{x}^{(\rho)}$, where $\rho=N$ or $\rho=N+n_0$, see Notation~\ref{Notation1}, and the respective Lagrange interpolation bases
 $\{\ell_{\rho-1, j}(x)\}_{j=1}^{\rho}$. Recall that  for each $j \in\{1, \ldots,\rho\}$, 
\begin{equation}
\ell_{\rho-1,j}(x)=\frac{\psi_\rho(x)}{\psi_\rho'\left(x_j^{(\rho)}\right)\left(x-x_j^{(\rho)}\right)}, 
\label{LagrangeBasisNnodes}
\end{equation}
where $\psi_\rho(x)=\left(x-x_1^{(\rho)}\right)\left(x-x_2^{(\rho)}\right)\cdots \left(x-x_\rho^{(\rho)}\right)$ is the node polynomial.

We define the generalized pseudospectral $(N+n_0) \times N$ matrix representation $A^{c}\left(\vec{x}^{(N)}, \vec{x}^{(N+n_0)}\right) \equiv A^c$ of the linear differential operator $\SA$  componentwise  by 
\begin{equation}
 A^{c}_{kj}= (\SA \ell_{N-1,j})\left(x_k^{(N+n0)}\right),  \;\; 1\leq k \leq N+n_0, 1\leq j \leq N,
\label{cRepres}
\end{equation}
 where, as before, the superscript `$c$' stands for ``collocation''. This definition is motivated by the following relation:
 $$
\pi_{N+n_0} \SA g= A^c \pi_N g
$$
for all $g \in \mathbb{P}^{N-1}$, where the isomorphisms $\pi_{\rho}: \mathbb{P}^{\rho-1} \to \mathbb{R}^{\rho}$ are defined by~\eqref{piN} with $N$ replaced by $\rho$, $\rho \in\{N, N+n_0\}$. 
In other words, if  ODE~\eqref{Au=f}  has a  polynomial solution $u \in \mathbb{P}^{N-1}$, then the vector $\vec{u}^{(N)} \in \mathbb{R}^N$ solves the system of linear equations 
$$A^c \vec{u}^{(N)}=\pi_{N+n_0}f$$
 if and only if the polynomial $u(x)=(\pi_N)^{-1} \vec{u}^{(N)}=\sum_{j=1}^N {u}^{(N)}_j \ell_{N-1, j}(x)$ solves ODE~\eqref{Au=f}.

Using analogous motivation, we define the generalized spectral $(N+n_0) \times N$  matrix representation $A^{\tau}$ of the linear differential operator $\SA$ componentwise by
\begin{equation}
A^{\tau}_{kj}=\frac{\langle \SA p_{j-1}, p_{k-1}\rangle}{\|p_{k-1}\|^2}, \;\; 1\leq k \leq N+n_0, 1\leq j \leq N.
\label{tauRepres}
\end{equation}
Here, the superscript `$\tau$' indicates that the $\tau$-variant of the spectral method is used~\cite{SpectralMethodsBook}.

We prove that the $(N+n_0) \times N$ matrices $A^\tau$ and $A^c$ satisfy the following property:
\begin{equation}
L^{(N+n_0-1)} A^c= A^\tau L^{(N-1)},
\label{AcAtauPseudoSimilarity}
\end{equation}
where each of the two $\rho \times \rho$ matrices $L^{(\rho-1)}$ with $\rho=N$  or $\rho=N+n_0$ is the transition matrix from the orthogonal polynomial basis $\{p_j(x)\}_{j=0}^{\rho-1}$ to the Lagrange interpolation basis $\{\ell_{\rho-1,j}\}_{j=1}^\rho$, see Theorem~\ref{ThmPseudoSimilarity} in Section~\ref{Sec:Proofs}.

\subsection{Main result: algebraic identities satisfied by the zeros of the polynomials $\{p_{\nu}(x)\}_{\nu=0}^\infty$}
\label{MainRes}

Let us now assume that $x_1^{(N)}, \ldots, x_{N}^{(N)}$ are the zeros of the polynomial $p_{N}(x)$ from the orthogonal family $\{p_\nu(x)\}_{\nu=0}^\infty$, while
 $x_1^{(N+n_0)}, \ldots, x_{N+n_0}^{(N+n_0)}$ are the zeros of the polynomial $p_{N+n_0}(x)$. We therefore use these two sets of zeros as the nodes  in the definition of the pseudospectral matrix representation $A^c$ of the differential operator $\SA$, see~\eqref{cRepres} and~\eqref{LagrangeBasisNnodes}.
In this case, each of the two matrices $L^{(\rho-1)}$ in the relation~\eqref{AcAtauPseudoSimilarity}, where $\rho=N$ or $\rho=N+n_0$, can be expressed in terms of the values $p_{m-1}(x_j^{(\rho)})$ and the Christoffel numbers $\lambda_j^{(\rho-1)}$:
\begin{eqnarray}
L^{(\rho-1)}_{mj}=\frac{p_{m-1}(x_j^{(\rho)})}{\|p_{m-1}\|^2} \lambda_j^{(\rho-1)}, \; 1\leq m,j \leq \rho,\label{LNm1pn0}
\end{eqnarray}
where the Christoffel numbers are defined by
\begin{equation}
\lambda_j^{(\rho-1)}=\int \ell_{\rho-1,j}(x) \, d\omega, \;\; 1\leq j\leq \rho,
\label{ChristoffelNumbers}
\end{equation}
 see Theorem~\ref{ThmPseudoSimilarity}. 

Recall that Christoffel numbers arise in the Gaussian quadrature numerical integration formulas; they are always positive~\cite{MastMilo}.  Christoffel numbers play an important role in the proof of the main identity~\eqref{eq1xnG} presented in this paper, although they are eliminated from that identity in the process of inversion of the matrices $L^{(\rho-1)}$: the inverses of $L^{(\rho-1)}$ are given componentwise by
\begin{equation}
\left\{\left[L^{(\rho-1)}\right]^{-1}\right\}_{mj}=p_{j-1}\left(x_m^{(\rho)}\right), \;  1\leq m,j \leq \rho.
\label{Lkm1Inverse}
\end{equation}
see~\eqref{ljtau} in the proof of Theorem~\ref{ThmPseudoSimilarity}.

Using property~\eqref{AcAtauPseudoSimilarity} of the matrix representations $A^\tau$ and $A^c$, together with the neat formulas~\eqref{Lkm1Inverse} for the matrices $\left[L^{(\rho-1)}\right]^{-1}$ in the case where $\vec{x}^{(\rho)}$ are the zeros of $p_\rho(x)$, $\rho \in \{N, N+n_0\}$, we prove the following algebraic identities satisfied by the zeros of the polynomials in the family $\{p_\nu(x)\}_{\nu=0}^\infty$.

\begin{theorem} 
The zeros $\vec{x}^{(N)}=\left(x_1^{(N)}, \ldots, x_{N}^{(N)}\right)$ of the polynomial $p_{N}(x)$ and the zeros $\vec{x}^{(N+n_0)}=\left(x_1^{(N+n_0)}, \ldots, x_{N+n_0}^{(N+n_0)}\right)$ of the polynomial $p_{N+n_0}(x)$ in the orthogonal polynomial family $\{p_{\nu}(x)\}_{\nu=0}^\infty$ of generalized eigenfunctions of the linear differential operator $\SA$, see~\eqref{EigenFunctionEquation}, satisfy the following algebraic relations for all integer $m,n$ such that  $1\leq m \leq N+n_0$, $1 \leq n \leq N$:
\begin{eqnarray}
&&\sum_{k=1}^{N} A^c_{mk}\left(\vec{x}^{(N)}, \vec{x}^{(N+n_0)}\right) p_{n-1}\left(x_k^{(N)}\right)\notag\\
&&=\sum_{ k=\max\{n-n_0,1\}}^{\min\{n+n_0,N+n_0\}} A^{\tau}_{kn} p_{k-1}\left( x_m^{(N+n_0)}\right),
\label{eq1xnG}
\end{eqnarray}
where the $(N+n_0)\times N$  pseudospectral and spectral matrix representations $A^c=A^c\left(\vec{x}^{(N)}, \vec{x}^{(N+n_0)}\right)\equiv A^c$ and  $A^\tau$, respectively,  are defined by~\eqref{cRepres} and~\eqref{tauRepres}.
\label{PropZerosThmG}
\end{theorem}

\begin{remark} For every pair of integers $n,m$ such that $0\leq n \leq N-1$ and $1\leq m \leq N+n_0$, identity~\eqref{eq1xnG} relates the zeros of the polynomials $p_N(x)$, $p_{N+n_0}(x)$ and $p_n(x)$ with the zeros of all the polynomials $p_k(x)$ such that the  index $k \in\{0,1, \ldots, N+n_0-1\}$ satisfies $n-n_0\leq k \leq n+n_0$.
\label{RemarkZerosNZerosmm1}
\end{remark}

\begin{remark} The main identity of Theorem~\ref{PropZerosThmG} may be recast as follows:
\begin{eqnarray}
  \left[\vec{e}^{(N+n_0)}_m\right]^T A^c \vec{v}^{(N)}_n=\vec{w}_m^{(N+n_0)} A^\tau \vec{e}_n^{(N)}
  \label{eq1xnGMatrixVecotorForm}
\end{eqnarray} 
for all $1\leq m \leq N+n_0$ and $1\leq n \leq N$,
where each $\vec{v}_n^{(N)} \in \mathbb{R}^N$ is a column-vector with the components $\left[\vec{v}_n^{(N)}\right]_k=p_{n-1}\left( x_k^{(N)}\right)$, each $\vec{w}_m^{(N+n_0)} \in \mathbb{R}^{N+n_0}$ is a row-vector with the components $\left[\vec{w}_m^{(N+n_0)}\right]_k=p_{k-1}\left( x_m^{(N+n_0)}\right)$ and the column-vectors $\{\vec{e}^{(\rho)}_k\}_{k=1}^\rho$ form the standard basis of $\mathbb{R}^{\rho}$.
\end{remark}

Theorem~\ref{PropZerosThmG} is proved in Section~\ref{Sec:Proofs}.

Identities~\eqref{eq1xnG} or, equivalently,~\eqref{eq1xnGMatrixVecotorForm}, are remarkable in the sense that they reveal a deeper structure and relation between the linear operators $A^c$ and $A^\tau$, the generalized spectral and pseudospectral representations, respectively, of the differential operator $\SA$, for the case where these operators are constructed using two sets of zeros of orthogonal polynomials as the nodes. 


Let us compare the main result of Theorem~\ref{PropZerosThmG} with other results of this kind. By setting $n_0=0$ in identities~\eqref{eq1xnG},~\eqref{eq1xnGMatrixVecotorForm}, we obtain that
\begin{equation}
\sum_{k=1}^N A^c_{mk}\left(\vec{x}^{(N)}\right) p_{n-1}\left(x_k^{(N)}\right)=q_{n-1} p_{n-1}(x_m^{(N)})
\label{IdentitieaKrall}
\end{equation}
for all $m,n$ such that $1\leq m,n \leq N$. In this case, of course, $q_\nu$ are constants, see~\eqref{EigenFunctionEquation}, and $\{p_\nu(x)\}_{\nu=0}^\infty$ are either classical or Krall orthogonal polynomials. In this special case where $n_0=0$,   identities of Theorem~\ref{PropZerosThmG} reduce to those reported in~\cite{Bihun2016}. It was shown in~\cite{Bihun2016} that, if applied to the classical Jacobi, Hermite or Laguerre polynomials,  identities~\eqref{IdentitieaKrall} reduce to known identities for the zeros of these polynomials reported in~\cite{Ahmed, Sasaki15}.

Therefore, identities~\eqref{eq1xnG},~\eqref{eq1xnGMatrixVecotorForm} generalize similar results for classical orthogonal polynomials proved in~\cite{Ahmed, Sasaki15} and for Krall polynomials proved  in~\cite{Bihun2016}. These identities may be considered as analogues of the properties of the zeros of the  Askey scheme and generalized hypergeometric polynomials  proved in~\cite{BCGenHyperg2014,BCAskey2014,BCGenBHyperg2015,BCqAskey2016}, for the case of the polynomial families considered in this paper.  An application of the identities proved in this paper is related to the study of the asymptotic behavior of algebraic expressions involving the zeros of orthogonal polynomials of degree $N$ as $N \to \infty$, see~\cite{Dominici2016,DominiciVanAssche2014}.

In the next Section~\ref{sec:SMpolys} we apply Theorem~\ref{PropZerosThmG} to prove new identities satisfied by the zeros of the nonclassical Sonin-Markov orthogonal polynomials. In Section~\ref{Sec:Proofs}, ``Proofs'', we elaborate on the proofs of most of  the theorems of this paper, except for those that are straightforward consequences of another theorem. In Section~\ref{Sec:Outlook} titled ``Discussion and Outlook'' we summarize the results proposed in this paper and discuss their importance, possible applications and further developments.

\section{Application: Properties of the Zeros of the Sonin-Markov Orthogonal Polynomials}
\label{sec:SMpolys}

In this Section~\ref{sec:SMpolys} we illustrate Theorem~\ref{PropZerosThmG} by applying it to the case of the  Sonin-Markov orthogonal polynomials, which are generalized eigenfunctions of a certain linear differential operator $\SD$, see~\eqref{diffEqSoninMarkov}. This application requires a computation of the generalized spectral and the generalized pseudospectral matrix representations of the differential operator $\SD$ associated with the Sonin-Markov polynomials. We thus proceed as follows. In Subsection~\ref{SM:Basics} we define the Sonin-Markov polynomials and state their basic properties. In Subsections~\ref{sec:DcSM} and~\ref{sec:DtauSM}, respectively, we compute the generalized pseudospectral and the generalized spectral representations, respectively, of the differential operator $\SD$. 
Finally,  in Subsection~\ref{SM:IdentitiesZeros} we provide a family of new algebraic properties satisfied by the zeros of the Sonin-Markov polynomials.

\subsection{Definition and Basic Properties of the Sonin-Markov Polynomials}
\label{SM:Basics}

The Sonin-Markov polynomials  $\{p^\beta_\nu(x)\}_{\nu=0}^\infty$, which are also known as generalized Hermite polynomials, are orthogonal on $I=(-\infty, \infty)$ with respect to the weight $w(x)=|x|^\beta e^{-x^2}$  (see~\cite{Kis,MastMilo} and references therein). These polynomials can be expressed in terms of the generalized Laguerre polynomials as follows:
\begin{subequations}
\begin{eqnarray}
p_{2n}^\beta(x)=c_n L_n^{\alpha}(x^2), \label{SMdfEven}\\
p_{2n+1}^\beta(x)=d_n x L_n^{\alpha+1}(x^2), \label{SMdfOdd}
\end{eqnarray}
\label{SMdf}
\end{subequations}
where 
\begin{equation}
\alpha=(\beta-1)/2
\label{alpha}
\end{equation} and the coefficients
\begin{eqnarray}
c_n=(-1)^n \sqrt{\frac{n!}{\Gamma(n+1+\alpha)}}, \; d_n=(-1)^n\sqrt{\frac{n!}{\Gamma(n+2+\alpha)}}
\label{cndnSM}
\end{eqnarray}
are chosen to make the polynomials orthonormal and the leading coefficients positive. Note that the leading coefficients $K_\nu$ of $p_{\nu}^\beta(x)$ are given by
\begin{eqnarray}
K_{2n}=\left[n! \Gamma(n+1+\alpha) \right]^{-1/2}, \notag\\
K_{2n+1}=\left[n! \Gamma(n+2+\alpha) \right]^{-1/2}.
\label{KnuSM}
\end{eqnarray}
The zeros of the Sonin-Markov polynomials are distinct, real and symmetric with respect to the origin.

The Sonin-Markov polynomials satisfy the following differential equations:
\begin{subequations}
\begin{eqnarray}
\mathcal{D} p_\nu(x)=(\mu_\nu x^2+\eta_\nu) p_\nu(x),
\end{eqnarray}
where $\mathcal{D}$ is the linear differential operator given by 
\begin{equation} 
\mathcal{D}=x^2 \frac{d^2}{dx^2}+x(-2x^2+2\alpha+1)\frac{d}{dx}
\label{DSoninMarkov}
\end{equation}
\label{diffEqSoninMarkov}
\end{subequations}
and
\begin{subequations}
\begin{eqnarray}
&&\mu_\nu=
-2\nu,\\
&&\eta_\nu=\left\{
\begin{array}{l}
0 \mbox{ if } \nu=2n,\label{etakSM}\\
2 \alpha+1 \mbox{ if } \nu=2n+1.
\end{array}
\right.
\end{eqnarray}
\label{muketakSM}
\end{subequations}
These differential equations can be derived from the corresponding differential equations satisfied by the generalized Laguerre polynomials stated in~\cite{KoekSwart}, see also~\cite{Szego} and formula (3.5) in~\cite{BCheikhGaied}.
Clearly, $\mathcal{D}\mathbb{P}^\nu \subseteq \mathbb{P}^{\nu+2}$, so we can apply Theorem~\ref{PropZerosThmG} with $n_0=2$ to these polynomials.

\subsection{The Generalized Pseudospectral Matrix Representation of the Sonin-Markov Differential Operator~\eqref{DSoninMarkov} }
\label{sec:DcSM}

Let  $x_1^{(N)},\ldots, x_N^{(N)}$ be the zeros of the Sonin-Markov polynomial $p_{N}(x) \equiv p_{N}^\beta(x)$ defined by~\eqref{SMdf} and let $x_1^{(N+2)},\ldots,  x_{N+2}^{(N+2)}$ be the zeros of the Sonin-Markov polynomial $p_{N+2}(x) \equiv p_{N+2}^\beta(x)$. In this subsection we find the generalized pseudospectral matrix representation $D^c\left( \vec{x}^{(N)}, \vec{x}^{(N+2)}\right)\equiv D^c$ of the Sonin-Markov differential operator~\eqref{DSoninMarkov} with respect to the two vectors of nodes $\vec{x}^{(N)}=\left(x_1^{(N)}, \ldots, x_N^{(N)}\right)$ and $\vec{x}^{(N+2)}=\left(x_1^{(N+2)}, \ldots, x_{N+2}^{(N+2)}\right)$. To compute this representation $D^c$, we use definition (\ref{cRepres}) with $n_0 = 2$. In the following, we use the notation of Subsection~\ref{SM:Basics}, note the appropriate definitions of $\alpha$, $c_n, d_n$ and $\mu_\nu, \eta_\nu$.

Let $\{\ell_{N-1,j}(x)\}_{j=1}^N$ be the Lagrange interpolation basis with respect to the nodes $\{x_j^{(N)}\}_{j=1}^N$ defined by (\ref{LagrangeBasisNnodes}) with  $\rho=N$ and $\psi_N(x)$ replaced by $p_N(x)$. Let  $Z^{(k)}\left( \vec{x}^{(N)}, \vec{x}^{(N+2)}\right)\equiv Z^{(k)}$ be the $(N+2)\times N$ generalized pseudospectral matrix representation of the differential operator $\frac{d^k}{dx^k}$ with respect to the two vectors of nodes $\vec{x}^{(N)}$ and $\vec{x}^{(N+2)}$. By definition~\eqref{cRepres}, its components $Z^{(k)}_{mn}$ are given by
\begin{equation}
Z^{(k)}_{mn}=\left[\frac{d^k}{dx^k} \ell_{N-1,n}(x) \right]\Bigg|_{x=x_m^{(N+2)}}.
\end{equation}

The generalized pseudospectral $(N+2)\times N$ matrix representation $D^c$ of the differential operator $\SD$ is given componentwise by
\begin{eqnarray}
D^c_{mn}&=&  \left[\SD\ell_{N-1,n}(x) \right]\Bigg|_{x=x_m^{(N+2)}} \notag\\
&=&\left[x_m^{(N+2)}\right]^2Z^{(2)}_{mn}+x_m^{(N+2)}\left\{-2\left[x_m^{(N+2)}\right]^2+2\alpha+1\right\} Z^{(1)}_{mn},
\label{DcmnNotSimplified}
\end{eqnarray}
see definition~\eqref{cRepres}.

By using the fact that the Sonin-Markov polynomial $p_{N}(x)$ satisfies differential equation~\eqref{diffEqSoninMarkov}  and the differential equation obtained by differentiation of~\eqref{diffEqSoninMarkov} with respect to $x$ (with $\nu=N$) , we simplify the formulas for $Z^{(k)}_{mn}$, where $k=1,2$, to read
\begin{equation}
Z^{(1)}_{mn}=
\begin{cases}
\frac{B_{mn}}{p'_N\big(x_n^{(N)}\big)}\Big[p'_N\big(x_m^{(N+2)}\big) - B_{mn}p_N\big(x_m^{(N+2)}\big)\Big]  \mbox{ if } ~ x_m^{(N+2)} \neq x_n^{(N)}, \\
0 \mbox{ if } ~ x_m^{(N+2)} = x_n^{(N)}=0
\end{cases}
\label{Z1}
\end{equation}
and
\begin{equation}
Z^{(2)}_{mn}=
\begin{cases}
\frac{B_{mn}}{p_N'\big(x_n^{(N)}\big)}\Bigg\{\left[\frac{2\left(x_m^{(N+2)}\right)^2 - 2\alpha - 1}{x_m^{(N+2)}} - 2B_{mn}\right]p'_N\big(x_m^{(N+2)}\big) \\ \qquad ~~~~~ 
+ \left[ \frac{\mu_N\left(x_m^{(N+2)}\right)^2 +\eta_N }{\left(x_m^{(N+2)}\right)^2}+ 2B_{mn}^2\right]p_N\big(x_m^{(N+2)}\big)\Bigg\} \\ \qquad ~~~~~ \mbox{if} ~ x_m^{(N+2)} \neq x_n^{(N)} \mbox{ and } x_m^{(N+2)} \neq 0, \\ 
\smallskip
\frac{-2p'_N(0)}{p'_N\big(x_n^{(N)}\big)}B_{mn}^2 \;\; \mbox{if} ~ x_m^{(N+2)} = 0 ~ \mbox{and } x_n^{(N)} \neq 0,\\ \smallskip
\frac{(\mu_N) + 2}{2(\alpha+2)} \;\; \mbox{if} ~ x_m^{(N+2)} = x_n^{(N)}=0,
\end{cases}
\label{Z2}
\end{equation}
where $B_{mn}$ are defined by 
\begin{equation}
\label{Amn}
B_{mn} = \left[x_m^{(N+2)} - x_n^{(N)}\right]^{-1}.
\end{equation}

By using the last two  expressions for the components of the matrices $Z^{(1)}$ and $Z^{(2)}$ in~\eqref{DcmnNotSimplified}, we obtain the following formulas for the components of the generalized pseudospectral matrix representation $D^c$:
\begin{subequations}
\label{Dc}
\begin{eqnarray}
D^c_{mn}& =& 
\frac{B_{mn}}{p'_N\big(x_n^{(N)}\big)}\Bigg[-2B_{mn}\big(x_m^{(N+2)}\big)^2p'_N\big(x_m^{(N+2)}\big) \notag \\ 
&&+  \bigg\{\mu_N\big(x_m^{(N+2)}\big)^2 + \eta_N + 2\big(B_{mn}\big)^2\big(x_m^{(N+2)}\big)^2 \notag\\ 
&& - B_{mn}x_m^{(N+2)}[-2\big(x_m^{(N+2)}\big)^2 + 2\alpha + 1] \bigg\}p_N\big(x_m^{(N+2)}\big)\Bigg] \notag\\ 
&&  \mbox{if} ~ x_m^{(N+2)} \neq x_n^{(N)} \mbox{ and } x_m^{(N+2)} \neq 0,  \label{Dca} \\
D^c_{mn} &=&  0  \;\; \;\;\;\; \mbox{if} ~ x_m^{(N+2)} = 0, \label{Dcb} 
\end{eqnarray}
\end{subequations}
where, again, $B_{mn}$ is given by definition (\ref{Amn}).

\begin{remark} If $\rho$ is even, the roots $\{x_j^{(\rho)}\}_{j=1}^\rho$ of the Sonin-Markov polynomial $p_\rho(x)\equiv p_\rho^\beta(x)$ are all distinct from zero because the roots of every generalized Laguerre polynomial are all distinct from zero, see definition~\eqref{SMdfEven}. Therefore, case~\eqref{Dcb} does not apply to even values of $N$. 
\end{remark}

\begin{remark} It is known that two generalized Laguerre polynomials $L_n^\alpha(x)$ and $L_{n+1}^\alpha(x)$ do not have common roots~~\cite{MastMilo,NikiforovUvarov}. But then, by definition~\eqref{SMdf}, two Sonin-Markov polynomials $p_{N}^\beta(x)$ and $p_{N+2}^\beta(x)$ have a common root $\hat{x}$ if and only if $N$ is odd and $\hat{x}=0$. This is why formulas~\eqref{Z1} and~\eqref{Z2}  include the case where $x_m^{(N+2)}=x_n^{(N)}=0$ (wich is possible for odd $N$ but not for even $N$), but do not include the case where $x_m^{(N+2)}=x_n^{(N)} \neq 0$ (the latter case is not possible).
\end{remark}

\subsection{The Generalized Spectral Matrix Representation of the Sonin-Markov Differential Operator~\eqref{DSoninMarkov} }
\label{sec:DtauSM}

The generalized spectral $(N+2)\times N$  matrix representation $D^\tau$ of the operator $\SD$ is defined by formula~\eqref{tauRepres}.  To find an explicit formula for the components of $D^\tau$, we will employ recurrence relations satisfied by the Sonin-Markov polynomials.

Recall that the generalized Laguerre polynomials satisfy certain three-term recurrence relations~\cite{MastMilo,NikiforovUvarov}, which imply 
\begin{eqnarray}
x^2 p_\nu(x)=\alpha_{\nu+2} p_{\nu+2}(x)+\beta_\nu p_\nu(x)+\gamma_{\nu-2} p_{\nu-2}(x)
\label{3termSM}
\end{eqnarray}
with
\begin{subequations}
\begin{eqnarray}
&&\alpha_{\nu}=\left\{
\begin{array}{l}
\left[n(n+\alpha)\right]^{1/2} \mbox{ if } \nu=2n,\\
\left[n(n+1+\alpha)\right]^{1/2}  \mbox{ if } \nu=2n+1,
\end{array}
\right. \label{alphakSM}\\
&&\beta_\nu=\nu+\alpha+1, \label{betakSM}\\
&&\gamma_{\nu}=\left\{
\begin{array}{l}
\left[(n+1)(n+1+\alpha)\right]^{1/2} \mbox{ if } \nu=2n,\\
\left[(n+1)(n+2+\alpha)\right]^{1/2}  \mbox{ if } \nu=2n+1,
\end{array}
\label{gammakSM}
\right.
\end{eqnarray}
\label{alphakbetakgammakSM}
\end{subequations}
where $p_k^\beta(x)$ are assumed to equal zero if $k<0$. Using recurrence relations~\eqref{3termSM}, we obtain
\begin{eqnarray}
&&D^{\tau}_{nj}=\eta_{n-1} \delta_{nj} +\mu_{j-1} 
\left(
\alpha_{n-1} \delta_{n, j+2} +\beta_{n-1} \delta_{n,j}+\gamma_{n-1} \delta_{n,j-2}
\right)
\end{eqnarray}
where $n=1,2, \ldots, N+2$, $j=1,2,\ldots, N$, and the coefficients  $\mu_\nu, \eta_\nu$ and $\alpha_\nu, \beta_\nu, \gamma_\nu$ are given by~\eqref{muketakSM} and~\eqref{alphakbetakgammakSM}. 

\begin{remark}
The coefficients $\alpha_\nu$ and $\beta_\nu$ defined in~\eqref{alphakSM} and~\eqref{betakSM}, respectively, must not be confused with the parameters $\alpha$ and $\beta$ in the definition of the Sonin-Markov polynomials, see~\eqref{SMdf} and~\eqref{alpha}.
\end{remark}

For completeness, let us mention that as a consequence of their orthogonality,  Sonin-Markov polynomials satisfy the standard three-term recurrence relation 
\begin{equation}
\label{srr}
xp_\nu(x) = \tilde{\alpha}_{\nu+1}p_{\nu+1}(x) + \tilde{\alpha}_{\nu}p_{\nu-1}(x) 
\end{equation}
where
\begin{equation}
\tilde{\alpha}_{\nu} = 
\begin{cases}
0 & \mbox{if} ~ \nu  = 0, \\
c_nd_n\frac{\Gamma(n+\alpha+2)}{n!} & \mbox{if} ~ \nu = 2n+1 ~ \mbox{and} ~ n \geq 0, \\
-c_nd_{n-1}\frac{\Gamma(n+\alpha+1)}{(n-1)!} & \mbox{if} ~ \nu = 2n ~ \mbox{and} ~ n \geq 1
\end{cases}
\end{equation}
and $p_j(x)=0$ if $j<0$. Of course, the three-term recurrence relation \eqref{3termSM} is a consequence of the recurrence relation \eqref{srr}.

\subsection{Algebraic Properties of the Zeros of the Sonin-Markov Polynomials}
\label{SM:IdentitiesZeros}

Having found the generalized pseudospectral and the generalized spectral matrix representations $D^c$ and $D^\tau$, respectively, of the differential operator $\SD$,  we apply Theorem~\ref{PropZerosThmG} to the case of the Sonin-Markov polynomials. We thus obtain the following algebraic identities satisfied by their zeros.
\begin{theorem}
For every pair of integers $n,m$ such that $1\leq n \leq N$ and $1 \leq m \leq N+2$, the zeros $x_1^{(N+2)}, x_2^{(N+2)},\ldots, x_{N+2}^{(N+2)}$ of the Sonin-Markov polynomial $p_{N+2}(x) \equiv p_{N+2}^{\beta}(x)$  and the zeros $x_1^{(N)},\ldots,x_N^{(N)}$ of the Sonin-Markov polynomial $p_{N}(x)\equiv p_{N}^{\beta}(x)$ satisfy the following relations. 

If $x_m^{(N+2)} = 0$, then
\begin{eqnarray}
&&\left(\eta_{n-1} + \mu_{n-1}\beta_{n-1}\right)p_{n-1}\left(0\right) + \mu_{n-1}\alpha_{n+1}p_{n+1}\left(0\right) \notag\\
&&+ \mu_{n-1}\gamma_{n-3}p_{n-3}\left(0\right) = 0,
\label{eq1xnSM1}
\end{eqnarray}
where $p_j(x)=0$ if $j<0$.

If $x_m^{(N+2)} \neq 0$, then
\begin{eqnarray}
&& \sum_{k=1}^{N}
\frac{B_{mk}\, p_{n-1}\left(x_k^{(N)}\right)}{p'_N\big(x_k^{(N)}\big)}\Bigg(-2B_{mk}\big(x_m^{(N+2)}\big)^2p'_N\big(x_m^{(N+2)}\big) \notag
\\ && +  \bigg\{\mu_N\big(x_m^{(N+2)}\big)^2 + \eta_N + 2\big(B_{mk}\big)^2\big(x_m^{(N+2)}\big)^2 \notag \\ && - \notag  B_{mk}x_m^{(N+2)}[-2\big(x_m^{(N+2)}\big)^2 + 2\alpha + 1] \bigg\}p_N\big(x_m^{(N+2)}\big)\Bigg)  \notag 
\\ && = \left(\eta_{n-1} + \mu_{n-1}\beta_{n-1}\right)p_{n-1}\left(x_m^{(N+2)}\right) \notag 
\\ && + \mu_{n-1}\alpha_{n+1}p_{n+1}\left(x_m^{(N+2)}\right) + \mu_{n-1}\gamma_{n-3}p_{n-3}\left(x_m^{(N+2)}\right) 
\label{eq1xnSM2}
\end{eqnarray}
where $B_{mn}$ are defined by~\eqref{Amn}
and, as before, $p_j(x)=0$ if $j<0$.
\label{PropZerosThmSM}
\end{theorem}

\section{Proofs}
\label{Sec:Proofs}
The proof of Theorem~\ref{PropZerosThmG} is based on the following result. 

\begin{theorem}
 Let $\mathcal{A}$ be a linear differential operator that satisfies condition~\eqref{DiffOperatorA}. Let $\{\ell_{N-1,j}(x)\}_{j=1}^N$ be the Lagrange interpolation basis of 
 $\mathbb{P}^{N-1}$ constructed using the $N$ distinct real nodes $x_1^{(N)}, \ldots, x_N^{(N)}$  and let $\{\ell_{N+n_0-1,j}(x)\}_{j=1}^{N+n_0}$ be the Lagrange interpolation basis of $\mathbb{P}^{N+n_0-1}$ constructed using the $(N+n_0)$ distinct real nodes $x_1^{(N+n_0)},  \ldots, x_{N+n_0}^{(N+n_0)}$, see~\eqref{LagrangeBasisNnodes}. 
If the $(N+n_0)\times N$ matrices $A^c, A^\tau$ are defined by~\eqref{cRepres} and~\eqref{tauRepres}, respectively, while the two $\rho \times \rho$ matrices $L^{(\rho-1)}$ with $\rho=N$ or $\rho=N+n_0$ are defined componentwise  by
\begin{eqnarray}
L^{(\rho-1)}_{mj}=\frac{\langle \ell_{\rho-1,j}, p_{m-1}\rangle}{\|p_{m-1}\|^2}, 1 \leq m,j \leq \rho,\label{LNPLambda}
\end{eqnarray}
then
\begin{equation}
L^{(N+n_0-1)} A^c=A^{\tau} L^{(N-1)}.
\label{Nice10}
\end{equation}
Moreover, if the interpolation nodes $x_1^{(\rho)}, \ldots, x_{\rho}^{(\rho)}$ are the distinct real zeros of the polynomial $p_{\rho}(x)$ from the orthogonal family $\{p_\nu(x)\}_{\nu=0}^\infty$ introduced in Section~\ref{Sec:Intro}, then the transition matrix $L^{(\rho-1)}$ and its inverse 
 $\left[L^{(\rho-1)}\right]^{-1}$ are given by~\eqref{LNm1pn0}  and~\eqref{Lkm1Inverse}, respectively, where $\rho \in\{N, N+n_0\}$.
 \label{ThmPseudoSimilarity}
\end{theorem}

\begin{remark} Let us note that $L^{(N-1)}$ is the transition matrix from the polynomial basis $\{p_m(x)\}_{m=0}^{N-1}$ to the basis $\{\ell_{N-1,m}(x)\}_{m=1}^N$ of $\mathbb{P}^{N-1}$, while $L^{(N+n_0-1)}$ is the transition matrix from the polynomial basis $\{p_m(x)\}_{m=0}^{N+n_0-1}$ to the basis $\{\ell_{N+n_0-1,m}(x)\}_{m=1}^{N+n_0}$ of $\mathbb{P}^{N+n_0-1}$.
\label{RemarkAppendix1}
\end{remark}

The proof of this theorem is similar to the proof of Theorem 1.1 in~\cite{Bihun2016}. It is provided below for the convenience of the reader.

\smallskip
\begin{proof}[Proof of Theorem \ref{ThmPseudoSimilarity}]

First, let us prove property~\eqref{Nice10}. Let $u$ be a polynomial of degree $N-1$. Then
\begin{eqnarray}
&&u(x)=\sum_{j=1}^N u_j^c \,\ell_{N-1,j}(x)\;\; \mbox{ and, on the other hand} \label{ucExpansion}\\
&&u(x)=\sum_{j=1}^N u_j^\tau \,p_{j-1}(x) \label{utauExpansion},
\end{eqnarray}
where the coefficients  
\begin{eqnarray}
&&u^c_j=u\left(x_j^{(N)}\right) \;\; \mbox{and} \label{uc}\\
&&u^\tau_j=\frac{\langle u, p_{j-1} \rangle}{\|p_{j-1}\|^2},
\label{utau}
\end{eqnarray}
respectively, are the components of the column-vectors $\vec{u}^c$ and $\vec{u}^\tau$, respectively, and $1\leq j \leq N$.
 To prove relation~\eqref{Nice10}, we will show that $\vec{u}^\tau = L^{(N-1)} \vec{u}^c$ and $L^{(N+n_0-1)} A^c \vec{u}^c= A^\tau \vec{u}^\tau$. 

Let us expand
\begin{equation}
\ell_{N-1,j}(x)=\sum_{m=1}^N L_{mj}^{(N-1)} p_{m-1}(x),
\label{ljtauExpansion}
\end{equation}
where the coefficients $L_{mj}^{(N-1)}$ are given by~\eqref{LNPLambda}. Upon a substitution of \eqref{ljtauExpansion}   into~\eqref{ucExpansion}, we obtain 
\begin{equation}
\vec{u}^\tau = L^{(N-1)} \vec{u}^c.
\label{utauuc}
\end{equation}

To obtain the equality $L^{(N+n_0-1)} A^c \vec{u}^c= A^\tau \vec{u}^\tau$, we first notice that because $u \in \mathbb{P}^{N-1}$ and the operator $\SA$ satisfies $\SA  \mathbb{P}^{N-1} \subseteq  \mathbb{P}^{N+n_0-1}$, we have
\begin{eqnarray}
\SA u(x)&=&\sum_{j=1}^{N+n_0} [A^c~\vec{u}^c]_j~\ell_{N+n_0-1,j}(x)=\sum_{j=1}^{N+n_0}  [A^c~\vec{u}^c]_j \sum_{m=1}^{N+n_0} L^{(N+n_0-1)}_{mj}~p_{m-1}(x)\notag\\
&=&\sum_{m=1}^{N+n_0} [L^{(N+n_0-1)} A^c~\vec{u}^c]_m~p_{m-1}(x).
\label{AucExpansion}
\end{eqnarray} 
On the other hand, 
\begin{eqnarray}
\SA u(x)&=&\SA \sum_{j=1}^N {u}^\tau_j p_{j-1}(x)=\sum_{j=1}^N  u^\tau_j \SA p_{j -1}(x)
=\sum_{j=1}^N  u^\tau_j \sum_{m=1}^{N+n_0} A^\tau_{mj} p_{m-1}(x)\notag\\
&=& \sum_{m=1}^{N+n_0} [A^\tau \vec{u}^\tau]_m p_{m-1}(x).
\label{AutauExpansion}
\end{eqnarray}
By comparing expansions~\eqref{AucExpansion} and~\eqref{AutauExpansion}, we obtain $L^{(N+n_0-1)} A^c \vec{u}^c= A^\tau \vec{u}^\tau$. Because $\vec{u}^\tau = L^{(N-1)}\vec{u}^c$, we conclude that $L^{(N+n_0-1)} A^c =A^\tau L^{(N-1)}$ and so finish the proof of relation~\eqref{Nice10}.

Second, let us assume that    $x_1^{(\rho)}, \ldots, x_\rho^{(\rho)}$ are the zeros of the polynomial $p_\rho(x)$, where $\rho=N$ or $\rho=N+n_0$. Let us prove that the transition matrix $L^{(\rho-1)}$ is given componentwise by~\eqref{LNm1pn0}. The Gaussian rule for approximate integration with respect to the measure $\omega$ based on these nodes $x_1^{(\rho)}, \ldots, x_\rho^{(\rho)}$ has degree of exactness $2\rho-1$, see, for example, Theorem 5.1.2 of~\cite{MastMilo}. Therefore, for  the polynomial $\ell_{\rho-1,j} (x) p_{m-1}(x)$ of degree $\rho-1+m-1\leq 2\rho-2$ we have
\begin{eqnarray}
&&\|p_{m-1}\|^2 L^{(\rho-1)}_{mj} = \int \ell_{N-1,j}(x) p_{m-1}(x) \, d\omega\notag \\
&&=\sum_{n=1}^\rho \ell_{N-1,j} (x_n^{(\rho)}) p_{m-1}(x_n^{(\rho)}) \lambda^{(\rho-1)}_n=p_{m-1}(x_j^{(\rho)}) \lambda^{(\rho-1)}_j,
\label{ljtau}
\end{eqnarray}
which implies~\eqref{LNm1pn0}. By applying the Gaussian rule to the polynomial $p_{m-1}(x)p_{n-1}(x)$, where $0\leq m,n \leq N$, we obtain
\begin{eqnarray}
\delta_{mn}=\int \frac{p_{m-1}(x) p_{n-1}(x)}{\|p_{m-1}\|^2} d\omega=\sum_{j=1}^\rho p_{n-1}(x_j^{(\rho)}) L_{mj}^{(\rho-1)},
\end{eqnarray}
which implies~\eqref{Lkm1Inverse}. 
\end{proof}

\begin{proof}[Proof of Theorem~\ref{PropZerosThmG}]
Theorem~\ref{PropZerosThmG} is a straightforward consequence of identity~\eqref{Nice10} rewritten as $\left\{A^c \left[L^{(N-1)}\right]^{-1}\right\}_{mn}=\left\{ \left[ L^{(N+n_0-1)} \right]^{-1} A^\tau \right\}_{mn}$, where $1\leq m \leq N+n_0$ and  $1\leq n \leq N$, and formula~\eqref{Lkm1Inverse} in Theorem~\ref{ThmPseudoSimilarity}.

We note that the index $k \in \{1,\ldots, N+n_0\}$ in the sum from the right hand side of identity~\eqref{eq1xnG} ranges from $n-n_0$ to $n+n_0$ rather than from $1$ to $N+n_0$. This is due to the following consideration. Because the polynomial family $\{p_\nu(x)\}_{\nu=0}^\infty$  is orthogonal and has the property that $\{p_{j}(x)\}_{j=0}^\nu$ is a basis of $\mathbb{P}^\nu$ for each $\nu$, this family must satisfy a three-term recurrence relation
\begin{equation}
x p_\nu(x)=a_{\nu, \nu+1} p_{\nu+1}(x)+a_{\nu, \nu} p_\nu(x) + a_{\nu, \nu-1} p_{\nu-1}(x),
\end{equation}
where $a_{\nu,k}$ are constants that equal zero if $k<0$ or $k$ is outside of the range $\{\nu-1, \nu, \nu+1\}$~\cite{MastMilo}.
From the last recurrence relation we derive
\begin{equation}
q_\nu(x) p_\nu(x)= \sum_{k=\nu-n_0}^{\nu+n_0} b_{\nu, k} p_{k}(x),
\end{equation}
where $b_{\nu, k}$ are constants that equal zero if $k<0$. Thus, for all integer $n,k$ such that $1\leq k \leq N+n_0$ and $1\leq n \leq N$,
$$A^\tau_{kn}=\frac{\langle q_{n-1} p_{n-1}, p_{k-1} \rangle}{\|p_{k-1}\|^2}=0 $$ 
if the index $k$ is outside of the range $\{n-n_0, \ldots, n+n_0\}$.
\end{proof}

\section{Discussion and Outlook}
\label{Sec:Outlook}

The identities of the main Theorem~\ref{PropZerosThmG} are derived using a very special \textit{exact} discretization of the differential equation 
\begin{equation}
\SA p_N(x)=q_N(x) p_N(x),
\label{DiffEqDisc1}
\end{equation}
compare with~\eqref{EigenFunctionEquation}. This discretization is constructed using the fact that ODE~\eqref{DiffEqDisc1} has a \textit{polynomial} solution $p_N(x)$, where $p_N(x)$ is a member of the nonclassical orthogonal polynomial family $\{p_\nu(x)\}_{\nu=0}^\infty$.  Because the differential operator $\SA$ maps $\mathbb{P}^N$ to $\mathbb{P}^{N+n_0}$, see property~\eqref{DiffOperatorA}, we choose to represent the operator by a $(N+n_0)\times N$ matrix $A^c\left(\vec{x}^{(N)}, \vec{x}^{(N+n_0)}\right)\equiv A^c$ defined by~\eqref{cRepres} in terms of \textit{two} vectors of interpolation nodes, $\vec{x}^{(N)}$ and $\vec{x}^{(N+n_0)}$. We thus \textit{generalize} the notion of pseudospectral matrix representations of a linear differential operator.

This generalization can be utilized to resolve the issues related to the incorporation of the initial or boundary conditions into pseudpsectral methods for solving a linear ODE $\SL u=f$, where $\SL$ is a linear differential operator, see~\cite{Bihun2011,BihunPrytula2010}. For example, if the differential operator $\SL$ has the property $\SL \mathbb{P}^\nu \subseteq  \mathbb{P}^{\nu-k_0}$, where $k_0$ is a positive integer, it is beneficial to use the generalized pseudospectral $N \times (N-k_0)$  matrix representation $L^c$ of the differential operator $\SL$ with respect to two vectors of nodes $\vec{x}^{(N)} \in \mathbb{R}^N$ and $\vec{x}^{(N-k_0)} \in \mathbb{R}^{N-k_0}$. Of course, $L^c$ is defined componentwise by $L^c_{kj}=(\SL \ell_{N-1,j})\left(x_k^{(N-k0)}\right),  \;\; 1\leq k \leq N-k_0, 1\leq j \leq N$, compare with~\eqref{cRepres}, where $\{\ell_{N-1, j}(x)\}_{j=1}^N$ is the standard Lagrange interpolation basis with respect to the nodes $\vec{x}^{N}$.

We can approximate the ODE $\SL u=f$ by the $N \times (N-k_0)$ system of linear equations $L^c \vec{u}^{(N)}=\vec{f}^{(N-k_0)}$, where $\vec{f}^{(N-k_0)}=\left(f\left(x_1^{(N-k_0)}\right),  \ldots, f\left(x_{N-k_0}^{(N-k_0)}\right) \right)$, compare with~\eqref{piN}. The linear system $L^c \vec{u}^{(N)}=\vec{f}^{(N-k_0)}$ for the $N$ components of $\vec{u}^{(N)}$ consists of $N-k_0$ linear equations, which allows to easily incorporate $k_0$ boundary or initial conditions into the system. Once the solution $\vec{u}^{(N)}$ of this augmented linear system is found, an approximate solution of the ODE $\SL u=f$ is given by $u(x)\approx (\pi_N)^{-1} \vec{u}^{(N)}=\sum_{j=1}^N {u}^{(N)}_j\ell_{N-1, j}(x)$, where $\{\ell_{N-1, j}(x)\}_{j=1}^N$ is the standard Lagrange interpolation basis with respect to the nodes $\vec{x}^{N}$. 

The ODE~\eqref{DiffEqDisc1} may also be discretized using only one vector of nodes $\vec{x}^{(N+n_0)}$ and the standard pseudospectral $(N+n_0)\times (N+n_0)$ matrix representation $\hat{A}^c$ of the differential operator $\SA$ defined componentwise by 
\begin{equation}
\hat{A}^c_{kj}=\SA \ell_{N+n_0-1, j}(x)\big|_{x=x_k^{(N+n_0)}}, \;\; 1\leq k,j \leq N+n_0,
\label{hatAc}
\end{equation}
see~\cite{Bihun2016} and compare with~\eqref{cRepres}. Indeed, in the proof of Theorem~\ref{ThmPseudoSimilarity} we may choose to represent the polynomial $u \in \mathbb{P}^{N-1}$ in terms of the Lagrange interpolation basis $\left\{\ell_{N+n_0-1,j}(x)\right\}_{j=1}^{N+n_0}$ or the basis $\{p_{j-1}(x)\}_{j=1}^{N+n_0}$ , see~\eqref{ucExpansion} and \eqref{utauExpansion}. This modification allows to prove the relation \begin{equation}
L^{(N+n_0-1)} \hat{A}^c= \hat{A}^\tau L^{(N+n_0-1)}, 
\label{RelationNpn0}
\end{equation}
where the standard spectral $(N+n_0) \times (N+n_0)$ matrix representation $\hat{A}^\tau$ is defined by 
\begin{equation}
\hat{A}^{\tau}_{kj}=\frac{\langle \SA p_{j-1}, p_{k-1}\rangle}{\|p_{k-1}\|^2}, \;\; 1\leq k,j \leq N+n_0.
\label{hatAtau}
\end{equation}
Let us assume again that $\vec{x}^{(N+n_0)}$ are the zeros of the polynomial $p_{N+n_0}$. Then relation~\eqref{RelationNpn0} rewritten as~$ \left\{\hat{A}^c\left[L^{(N+n_0-1)}\right]^{-1}\right\}_{mn}= \left\{\left[L^{(N+n_0-1)}\right]^{-1}\hat{A}^\tau \right\}_{mn}$ together with expression~\eqref{Lkm1Inverse} for the inverse of the transition matrices $L^{(N+n_0-1)}$ yield the identities
\begin{equation}
\sum_{k=1}^{N+n_0} \hat{A}^c_{mk}(\vec{x}^{(N+n_0)}) p_{n-1}\left(x_k^{(N+n_0)} \right)=\sum_{j=\max\{1,n-n_0\}}^{\min\{N+n_0, n+n_0\}} \hat{A}^{\tau}_{jn} p_{j-1}\left(x_m^{(N+n_0)} \right),
\label{Lalala}
\end{equation}
which hold for all integer $m,n$ such that $1\leq m,n \leq N+n_0$. Given a fixed integer $n$ such that $0\leq n \leq N+n_0-1$, the last identity relates the zeros of the polynomials $p_{N+n_0}(x)$ and $p_n(x)$ with the zeros of all the polynomials $p_j(x)$ such that the integer index $j \in \{0,1,\ldots, N+n_0-1\}$ satisfies $n-n_0 \leq j \leq n+n_0$.

The set of identities~\eqref{Lalala}  is a straighforward generalization of Theorem~1.1 in~\cite{Bihun2016}, thus it is not the main focus of this paper. We invite the reader to apply  identities~\eqref{Lalala} to the Sonin-Markov polynomials, replacing $ \hat{A}^c$ and  $\hat{A}^{\tau}$ with  $\hat{D}^c$ and $\hat{D}^\tau$, respectively, the latter being matrix representations of the differential operator~$\SD$ associated with the Sonin-Markov family, see~\eqref{DSoninMarkov}. The pseudospectral matrix representation $\hat{D}^c$ can be found using definition~\eqref{hatAc} and the techniques outlined in Appendix~A of~\cite{Bihun2016}. The spectral matrix representation $\hat{D}^\tau$ can be found using definition~\eqref{hatAtau} and recurrence relations~\eqref{3termSM}. 

In summary, the method of generalized pseudospectral matrix represenetations proposed in this paper can be utilized to prove new identities satisfied by nonclassical orthogonal polynomials that are generalized, rather than standard, eigenfunctions of linear differential operators. The proposed generalized matrix representations of linear differential operators  are quite useful in the process of incorporation of initial or boundary conditions into a discretization of a linear ODE via a pseudospectral method. The last application deserves further exploration involving numerical tests, which the authors may pursue in the near future. Other interesting developments of the method proposed in this paper may involve its generalization to exceptional orthogonal polynomials~\cite{GomezUlateCo2014, GomezUlateCo2015, OS2011}.

\section*{Acknowledgements}

This work was supported by a grant awarded to O. Bihun by the Committee on Research and Creative Works (CRCW) of the University of Colorado, Colorado Springs.

The authors declare that there is no conflict of interest regarding the publication of this paper.

\end{document}